\documentclass[11pt]{article}
\usepackage{amsfonts}
\usepackage{amssymb}

\newtheorem{theorem}{Theorem}[section]

\newtheorem{lemma}[theorem]{Lemma}
\newtheorem{corollary}[theorem]{Corollary}

 \usepackage{amssymb}
\usepackage{epsfig}
\usepackage{amsmath}
\usepackage{amsfonts}

\newcommand{\R}{{{\Bbb R}}}

\def\qed{\hbox to 0pt{}\hfill$\rlap{$\sqcap$}\sqcup$\medbreak}

\title{Existence and computation of Riemann--Stieltjes integrals through Riemann integrals}
\date{July, 2011}
\begin{document}
\maketitle

\vspace{1cm}

\begin{center}
{\large Rodrigo L\'opez Pouso \\
Departamento de
An\'alise Matem\'atica\\
Facultade de Matem\'aticas,\\Universidade de Santiago de Compostela, Campus Sur\\
15782 Santiago de
Compostela, Spain.
}
\end{center}

\begin{abstract}
We study the existence of Riemann--Stieltjes integrals of bounded functions against a given integrator. We are also concerned with the possibility of computing the resulting integrals by means of related Riemann integrals. In particular, we present a new generalization to the well--known formula for continuous integrands and continuously diffe\-ren\-tiable integrators.
\end{abstract}

 \section{Introduction}
Let $f:[a,b]\longrightarrow  \R$ be continuous and let $G:[a,b]\longrightarrow \R$ have bounded variation over the interval $[a,b]$. A standard result \cite[Exercise 30 (g), page 281]{str} then guarantees that
\begin{equation}
\label{rs}
\int_a^b{f(x) \, dG(x)}
\end{equation}
exists in the Riemann--Stieltjes sense. Moreover, the formula
\begin{equation}
\label{for1}
\int_a^b{f(x) \, dG(x)}=\int_a^b{f(x)G'(x) \, dx}
\end{equation} 
holds true if $G$ has a continuous derivative in $[a,b]$, thus reducing the computation of Riemann--Stieltjes integrals to Riemann ones.

\bigbreak 

We can find deeper information in Stromberg \cite[Exercise 30 (k), page 281]{str}: the formula (\ref{for1}) is fulfilled when $f$ is Riemann integrable and $G$ is absolutely continuous. This result leans on the Fundamental Theorem of Calculus for the Lebesgue integral, and the proof suggested by Stromberg uses approximations of $f$ by step functions. Notice that, in this case, the integral in the right--hand side of (\ref{for1}) is a Lebesgue integral which need not be, in general, a Riemann integral.

\medbreak

In this note we are concerned with the twofold problem of the existence of (\ref{rs}) in the Riemann--Stieltjes sense and the possibility of computing it via (\ref{for1}) when $f$ is merely a bounded function (not necessarily Lebesgue measurable). Obviously, the conditions over the integrator $G$ have to be reenforced. Specifically, in this paper we prove the following result.

\begin{theorem}
\label{th}
Let $f:[a,b]\longrightarrow \R$ be bounded, let $g:[a,b]\longrightarrow \R$ be Riemann integrable and let $G(x)=c+\int_a^x{g(y) \, dy}$ ($x \in [a,b]$) for some $c \in \R$.

A neccessary and sufficient condition for the function $f$ to be Riemann--Stieltjes integrable with respect to $G$ over $[a,b]$ is that the product $f g$ be Riemann integrable over $[a,b]$ and, in that case,
\begin{equation}
\label{rsf}
\int_a^b{f(x) \, dG(x)}=\int_a^b{f(x) g(x) \, dx}.
\end{equation}

\end{theorem}
 The proof of Theorem \ref{th}, which occupies section 2, is based on the following sharp version of the mean value theorem for Riemann integrals.

\begin{theorem}
\label{corpro}{\bf \cite[Corollary 4.6]{rlp}}
If $h:[a,b] \longrightarrow \R$ is Riemann integrable on $[a,b]$ then there exist points $c_1,\, c_2 \in (a,b)$ such that
\begin{equation}
\nonumber
h(c_1)(b-a) \le  {\int_a^b}{h(x) \, dx}  \le h(c_2)(b-a).
\end{equation}
\end{theorem}

Even though our conditions on the function $f$ are very general, our proofs turn out to be quite easy (this fact being another interesting point). In particular, we remain in the realm of Riemann integration theory.

\section{Proof of Theorem \ref{th}}
We need the following lemma on mixed Riemann sums (see \cite{cat} for more information on mixed sums).
Its proof is very easy, but we include it for completeness and for the convenience of the reader.

\begin{lemma}
\label{mixed}
Let $f:[a,b]\longrightarrow \R$ be bounded and let $g:[a,b]\longrightarrow \R$ be Riemann integrable.

If the product $f g$ is Riemann integrable on $[a,b]$ then for each $\varepsilon>0$ there exists $\delta>0$ such that for every partition $a=x_0<x_1<\dots<x_n=b$ whose norm is less than $\delta$ we have
$$\left|\sum_{k=1}^n{f(y_k)g(z_k)(x_k-x_{k-1})-\int_a^b{f(x)g(x) \, dx}} \right|<\varepsilon,$$
for any choice of points $y_k,z_k \in [x_{k-1},x_k]$ ($k=1,2,\dots,n$).
\end{lemma}

\noindent
{\bf Proof.} For a given partition $a=x_0<x_1<\dots<x_n=b$ and points $y_k,z_k \in [x_{k-1},x_k]$ ($k=1,2,\dots,n$) we have
\begin{align*}
\left|\sum_{k=1}^nf(y_k)g(z_k)\right.&\left.(x_k-x_{k-1})-\int_a^b{f(x)g(x) \, dx} \right| \\
&\le \left|\sum_{k=1}^n{f(y_k)[g(z_k)-g(y_k)](x_k-x_{k-1})}\right|\\
& \qquad +\left|\sum_{k=1}^n{f(y_k)g(y_k)(x_k-x_{k-1})} -\int_a^b{f(x)g(x) \, dx} \right| \\
&\le \sup_{a \le x \le b}|f(x)|\sum_{k=1}^n{\mbox{osc}(g,[x_{k-1},x_k])(x_k-x_{k-1})}\\
&\qquad+\left| \sum_{k=1}^n{f(y_k)g(y_k)(x_k-x_{k-1})}-\int_a^b{f(x)g(x) \, dx} \right|,
\end{align*}
where $$\mbox{osc}(g,[x_{k-1},x_k])=\sup_{x_{k-1}\le x \le x_k} g(x) -\inf_{x_{k-1}\le x \le x_k} g(x)$$ is the oscillation of $g$ in $[x_{k-1},x_k]$.

Since $g$ and $fg$ are Riemann integrable, the last term in the previous inequality is as small as we wish if the norm of the partition is sufficiently small (and this does not depend on the choice of $y_k, z_k \in [x_{k-1},x_k]$).\qed

\bigbreak

\noindent
{\bf Proof of Theorem \ref{th}.} We first show that the condition is sufficient, so we assume that $fg$ is Riemann integrable on $[a,b]$. Let $\varepsilon>0$ be fixed and let $\delta>0$ be as in Lemma \ref{mixed}. We are going to prove that for any partition of $[a,b]$, say $P=\{x_0,x_1,\dots,x_n\}$,  with norm less than $\delta$, and any choice of points $y_k \in [x_{k-1},x_k]$ ($k=1,2,\dots,n$), we have
 \begin{equation}
 \label{eq1}
 \left|\sum_{k=1}^n{f(y_k)[G(x_k)-G(x_{k-1})]}-\int_a^b{f(x)g(x) \, dx} \right|< \varepsilon,
 \end{equation}
 thus finishing the first part of the proof of Theorem \ref{th}.
 
 Let $P$ and the $y_k$'s be as above. Theorem \ref{corpro} guarantees that for each $k \in \{1,2,\dots,n\}$ there is some $z_k \in (x_{k-1},x_k)$ such that
 $$f(y_k) \, \int_{x_{k-1}}^{x_k}{g(y) \, dy} \le f(y_k)g(z_k)(x_k-x_{k-1}).$$
 Hence
\begin{align*}
\sum_{k=1}^n{f(y_k)[G(x_{k})-G(x_{k-1})]}&=\sum_{k=1}^n{f(y_k)\int_{x_{k-1}}^{x_k}{g(y) \, dy}} \\
&\le \sum_{k=1}^n{f(y_k)g(z_k)(x_k-x_{k-1})}\\
&< \int_a^b{f(x)g(x)\, dx}+\varepsilon \quad \mbox{(by Lemma \ref{mixed}).}
\end{align*}
We deduce in an analogous way that
$$\sum_{k=1}^n{f(y_k)[G(x_{k})-G(x_{k-1})]}> \int_a^b{f(x)g(x)\, dx}-\varepsilon,$$
so (\ref{eq1}) obtains.

\bigbreak

We now prove that the condition is necessary, and therefore we assume that $f$ is Riemann--Stieltjes integrable with respect to $G$ over $[a,b]$. Let us consider an arbitrary partition $a=x_0<x_1<\dots<x_n=b$ and arbitrary points $y_k \in [x_{k-1},x_k]$ ($k=1,2,\dots,n$). We have
\begin{align} 
\label{in1}
&\left| \sum_{k=1}^n { f(y_k)g(y_k)(x_k-x_{k-1}) } - \sum_{k=1}^n{f(y_k)[G(x_{k})-G(x_{k-1})]} \right| \\
\nonumber
& \qquad=\left|\sum_{k=1}^n {f(y_k)\left[g(y_k)(x_k-x_{k-1})- \int_{x_{k-1}}^{x_k}{g(y) \, dy} \right]}\right| \\
\label{in2}
& \qquad \le \sup_{a \le x \le b}|f(x)| \sum_{k=1}^n{\left| g(y_k)(x_k-x_{k-1})  -\int_{x_{k-1}}^{x_k}{g(y) \, dy}  \right|.      }
\end{align}
 For each $k \in \{1,2,\dots,n\}$ we have 
$$ \inf_{x_{k-1}\le x \le x_k}g(x) \, (x_k-x_{k-1}) \le \int_{x_{k-1}}^{x_k}{g(y) \, dy} \le \sup_{x_{k-1}\le x \le  x_k}g(x) \, (x_k-x_{k-1}),$$
and, since $y_k \in [x_{k-1},x_k]$, we deduce that
\begin{align*}
&\left| g(y_k)(x_k-x_{k-1})  -\int_{x_{k-1}}^{x_k}{g(y) \, dy}  \right| \\
& \qquad \le \left[ \sup_{x_{k-1}\le x_k}g(x)-\inf_{x_{k-1}\le x \le x_k}g(x) \right] (x_k-x_{k-1})\\
&  \qquad =\mbox{osc}(g,[x_{k-1},x_k])(x_k-x_{k-1}).
\end{align*}
Going with this information back to the inequality (\ref{in1})--(\ref{in2}), we get
\begin{align*}
&\left| \sum_{k=1}^n { f(y_k)g(y_k)(x_k-x_{k-1}) } - \sum_{k=1}^n{f(y_k)[G(x_{k})-G(x_{k-1})]} \right| \\
& \quad \le \sup_{a \le x \le b}|f(x)| \sum_{k=1}^n{\mbox{osc}(g,[x_{k-1},x_k])(x_k-x_{k-1})},
\end{align*}
which tends to zero when the norm of the partition tends to zero (and this does not depend on the $y_k$'s).

To sum up, if the norm of the partition is sufficiently small then Riemann sums of $fg$ are as close as we wish to Riemann--Stieltjes sums of $f$ with respect to $G$, which, in turn, are as close as we wish to the corresponding Riemann--Stieltjes integral. Hence Riemann sums of $fg$ are as close as we wish to the Riemann--Stieltjes integral of $f$ with respect to $G$ provided that the norm of the partition is sufficiently small.
\qed
 \section{Concluding remarks}
 Theorem \ref{th} is not a particular case to \cite[Exercise 30 (k)]{str}, which requires $f$ to be Riemann integrable on $[a,b]$. Examples in the conditions of Theorem \ref{th} which do not fulfill the assumptions in \cite[Exercise 30 (k)]{str} can be easily constructed: it suffices to consider bounded functions $f$ which are not Riemann integrable on a subinterval of $[a,b]$ where $g$ is almost everywhere equal to zero.

 If we combine Theorem \ref{th} with the {\it integration by parts formula} for the Riemann--Stieltjes integral, see \cite[Exercise 30 (h)]{str}, we obtain the following new sufficient condition for Riemann--Stieltjes integrability.
 
 \begin{corollary}
 \label{coth}
 In the conditions of Theorem \ref{th}, the function $G$ is Riemann--Stieltjes integrable with respect to $f$ over $[a,b]$ and
 $$\int_a^b{G(x) \, df(x)}=G(b)f(b)-G(a)f(a)-\int_a^b{g(x)f(x) \, dx}.$$ 
 \end{corollary} 
 
 Finally, Theorem \ref{th} and Corollary \ref{coth} yield the following more ``symmetric" test.
 
 \begin{corollary}
 Let $\alpha, \, \beta:[a,b] \longrightarrow \R$ be Riemann integrable functions.
 
 If one of them is an indefinite Riemann integral, then either of them is Riemann--Stieltjes integrable with respect to the other over $[a,b]$ and the corresponding Riemann--Stieltjes integrals reduce to Riemann integrals.
  \end{corollary}

\end{document}